\documentclass[11pt]{article}
\usepackage{amsmath,amssymb,amsthm,geometry,verbatim,enumerate,float,tikz}
\newtheorem{thm}{Theorem}
\newtheorem{lemma}[thm]{Lemma}
\newtheorem{coro}[thm]{Corollary}

\usepackage{setspace}

\begin{document}

%\linenumbers
\onehalfspace

\title{Random Subgraphs in Sparse Graphs}

\author{Felix Joos}

\date{}

\maketitle

\begin{center}
Institut f\"{u}r Optimierung und Operations Research, 
Universit\"{a}t Ulm, Ulm, Germany\\
\texttt{felix.joos@uni-ulm.de}
\end{center}

\begin{abstract}
We investigate the threshold probability for connectivity of sparse graphs under weak assumptions.
As a corollary this completely solve the problem for Cartesian powers of arbitrary graphs.
In detail, let $G$ be a connected graph on $k$ vertices, 
$G^n$ the $n$-th Cartesian power of $G$,
$\alpha_i$ be the number of vertices of degree $i$ of $G$, 
$\lambda$ be a positive real number, 
and $G^n_p$ be the graph obtained from $G^n$
by deleting every edge independently with probability $1-p$.
If $\sum_{i}\alpha_i(1-p)^i=\lambda^{\frac{1}{n}}$,
then $\lim_{n\rightarrow \infty}\mathbb{P}[G^n_p {\rm\ is\ connected}]=\exp(-\lambda)$.
This result extends known results for regular graphs.
The main result implies that the threshold probability does not depend on the graph structure of $G$ itself,
but only on the degree sequence of the graph.

\bigskip

%{\bf AMS subject classification:} 
%05C38, % Paths and cycles
%05C55 % Generalized Ramsey theory
\end{abstract}

\section{Introduction}
For a graph $G$ and $p\in(0,1)$ let $G_p$ be the graph obtained from $G$ by deleting  every edge independently
with probability $1-p$. 
First introduced by Erd{\H{o}}s and R\'enyi \cite{erdosrenyi}, such so-called random graphs are studied in great detail. 
For a survey see for example \cite{bolobook}.
Erd{\H{o}}s and R\'enyi showed that the probability that $(K_n)_p$ is connected
tends to $\exp(-e^{-c})$ if $p=\frac{\ln n +c}{n}$ for $c\in\mathbb{R}$ as $n\rightarrow \infty$.
There are similar results for complete bipartite graphs $K_{n,n}$ \cite{palasti}
and even for multipartite graphs \cite{ruci}. 

For a graph $G$, the $n$-th Cartesian power $G^n$ has vertex set $V(G)\times \ldots \times V(G)$
and two vertices $(v_1,\ldots, v_n)$ and $(w_1,\ldots, w_n)$
are adjacent if and only if there is an $i\in[n]$
such that $v_j=w_j$ for all $j\not= i$ and $v_iw_i\in E(G)$.
The best known Cartesian power of a graph is the $n$-dimensional hypercube $K_2^n$.
Burtin \cite{burtin}, Erd{\H{o}}s and Spencer \cite{erdosspencer}, and Bollob\'as \cite{bolo} showed that 
$\lim_{n\rightarrow \infty}\mathbb{P}[(K_2^n)_p {\rm\ is\ connected}]=\exp(-e^{-c})$
if $p=\frac{1}{2}+\frac{c}{2n}$ for every $c\in \mathbb{R}$.
Recently van der Hofstad and Nachmias \cite{hofstad} investigated the behaviour of the giant component in the hypercube.
Clark \cite{clark} considered Cartesian powers of complete and complete bipartite graphs.
He showed that $\lim_{n\rightarrow \infty}\mathbb{P}[(K_k^n)_p {\rm\ is\ connected}]=\exp(- \lambda)$
if $k\geq 2$, $1-p= \left(\frac{(\lambda_n)^{1/n}}{k}\right)^{\frac{1}{k-1}}$, and $\lambda_n\rightarrow \lambda>0$ as $n\rightarrow \infty$.
In addition, he proved an analogous result for $K_{k,k}$ if $k\geq1$.
Joos \cite{jo} generalized these results by proving that
$\lim_{n\rightarrow \infty}\mathbb{P}[G^n_p {\rm\ is\ connected}]=\exp(-\lambda)$
if $G$ is a connected $d$-regular graph ($d\geq 1$) on $k$ vertices,
$1-p= \left(\frac{(\lambda_n)^{1/n}}{k}\right)^{\frac{1}{d}}$, and $\lambda_n\rightarrow \lambda>0$ as $n\rightarrow \infty$.
Up to now nothing was known about the case if $G$ is non-regular - even if $G=P_3$.
In this paper we solve the problem for all graphs.

However, we are interested in a much general setting.
Instead of considering the sequence $G^1, G^2, \ldots$ of all Cartesian powers of some graph $G$,
we consider sequences of graphs $G_1, G_2, \ldots$ such that $G_n$ has much less structure than the $n$-th Cartesian power of some graph.
Let the random variable $X_n$ be the number of isolated vertices in $(G_n)_p$.
We prove, if $p$ is chosen such that $\mathbb{E}[X_n]=\lambda_n$ and $\lambda_n\rightarrow \lambda>0$ as $n\rightarrow \infty$,
then $\lim_{n\rightarrow \infty}\mathbb{P}[(G_n)_p {\rm\ is\ connected}]=\exp(-\lambda)$.

For a graph $H$, let $\delta(H)$ and $\Delta(H)$ be the minimum and maximum degree of $H$, respectively.
For a set of vertices $S$ in $H$, let $b_H(S)$ be number of edges that join $S$ and its complement in $H$
and call this set of edges the \textit{boundary} of $S$.
Let $b_H(s)=\min_{|S|=s}b_H(S)$.

We say a sequence $G=(G_n)_{n\geq 1}$ of connected graphs satisfies the \textit{basic} conditions for some $k\geq 2$, if
\begin{enumerate}
	\item $n(G^n)=k^n$,
	\item $\delta(G_n)\geq n$,
	\item $\exists \epsilon'=\epsilon'(G)>0: \ \frac{\Delta(G^n)}{\delta(G^n)}\leq n^{1-\epsilon'}$,
	\item $\exists c=c(G)\in\mathbb{N}:\ \Delta(G^n)\leq n^{c}$ and
	\item $\exists \epsilon=\epsilon(G)>0\ \forall s\in \{n^{1- \epsilon'},\ldots,k^n/2\}: \ b_{G^n}(s)\geq \epsilon \Delta(G^n) s\left(1- \frac{1}{n}\log_k s\right)$.
\end{enumerate}
Note that condition 5 is very natural for sparse graphs and
weaker than a condition of the type $C\Delta(G^n)s\left(1- \frac{s}{k^n}\right)$ for some constant $C$.
Our main contribution is Theorem \ref{mainresult1}.
\begin{thm}\label{mainresult1}
Let $G=(G_n)_{n\geq 1}$ be a sequence of connected graphs that satisfies the basic conditions for some
$k\in \mathbb{N}\setminus \{1\}$ and let $\lambda_n\rightarrow \lambda>0$ as $n\rightarrow\infty$.
If $p=p(n)$ is chosen such that $\mathbb{E}[X_n]=\lambda_n$,
then
\begin{align*}
	\lim_{n\rightarrow \infty}\mathbb{P}[(G_n)_p {\rm\ is\ connected}]=\exp(-\lambda).
\end{align*}
\end{thm}

For a connected graph $H$, let $P(H,x)$ be the polynomial $\sum_i\alpha_{i,H}x^i$,
where $\alpha_{i,H}$ is the number of vertices of degree $i$ in $H$.
Since all coefficients are non-negative, the equation $P(H,x)=y$ has a unique positive solution if $y>0$.

\begin{coro}\label{main result}
Let $H$ be a connected graph on $k\geq 2$ vertices. 
Let $\lambda_n\rightarrow \lambda>0$ as $n\rightarrow\infty$ 
and let $q$ be the unique positive solution of the equation $P(H,x)=\lambda^{\frac{1}{n}}_n$.
If $p=1-q$,
then 
\begin{align*}
	\lim_{n\rightarrow \infty}\mathbb{P}[H^n_p {\rm\ is\ connected}]=\exp(-\lambda).
\end{align*}
\end{coro}

\noindent
Note that, if $H$ is disconnected, then $H^n$ is also disconnected.
Furthermore, Corollary \ref{main result} implies that the
threshold function for connectivity does not depend on the graph structure,
but only on the degree sequence of $H$.

If $H$ is a $d$-regular graph,
then $P(H,x)=kx^d$ and $1-p=\left(\frac{\lambda_n^{1/n}}{k}\right)^{\frac{1}{d}}$.
Thus Corollary \ref{main result} extends all mentioned former results concerning Cartesian powers.
If $P_k$ is the path on $k$ vertices,
then $P(P_k,x)=(k-2)x^2+2x$ and hence $1-p=\frac{1}{k-2}\left(-1+\sqrt{1+(k-2)\lambda_n^{1/n}}\right)$.

%%%%%%%%%%%%%%%
\section{Preliminaries}
%%%%%%%%%%%%%%%

For a graph $G$, we denote by $V(G)$ and $E(G)$ its vertex and edge set, respectively.
Let the order $n(G)$ of $G$ be the number of vertices of $G$
and the size $m(G)$ of $G$ be the number of edges of $G$.
For a vertex $v\in V(G)$, let the neighborhood $N_G(v)$ be the set of all neighbors of $v$ in $G$
and the closed neighborhood $N_G[v]$ be defined by $N_G(v)\cup \{v\}$.
Let the degree $d_G(v)$ of $v$ in $G$ be defined by $|N_G(v)|$.
We denote by $\delta(G)$ and $\Delta(G)$ the minimum and maximum degree of $G$, respectively.
Let $S\subseteq V(G)$.
We denote by $G[S]$ the subgraph of $G$ induced by $S$.
Let $N_G(S)$ be the set of vertices in $V(G)\setminus S$ with at least one neighbor in $S$ and 
let $N_G[S]=N_G(S)\cup S$.
We say that $S$ dominates $G$ if $N_G[S]=V(G)$.
%Let the domination number $\gamma(G)$ of $G$ be the smallest number of vertices in a set $S$ such that $N_G[S]=V(G)$.
For a pair of vertices $u,v\in V(G)$, let the smallest number of edges on a $u,v$-path be the distance between $u$ and $v$.
For $\ell\in \mathbb{N}$,
a set $S$ of vertices $\ell$-dominates $G$
if every vertex $v$ of $G$ is either contained in $S$
or there is a vertex in $S$ in distance at most $\ell$ to $v$.
Let the $\ell$-domination number $\gamma^{\ell}(G)$ of $G$ be the smallest number of vertices in a set $S$ such that $S$ $\ell$-dominates $G$.
Note that $1$-domination coincides with domination defined above.
For a set $X$, let $\binom{X}{\ell}$ be the set of subsets of $X$ of cardinality $\ell$.

%Let the size of the boundary $b_G(S)$ of a set $S\subseteq V(G)$ be the number of edges
%that join $S$ and $V(G)\setminus S$.
%For $s\in \mathbb{N}$, let $b_G(s)$ be smallest size of a boundary of a set $S\subseteq V(G)$ such that $|S|=s$.
%If $(v_1,\ldots,v_n)\in V(G^n)$,
%then $d_{G^n}((v_1,\ldots,v_n))=d_G(v_1)+\ldots+d_G(v_n)\geq \delta(G) n$.

%We write $\ln x$ instead of $\log_e x$.
%Let $X$ be a random variable.
%We denote by $\mathbb{E}[X]$ the expected value of $X$ and let
%the $r$-th factorial moment $\mathbb{E}_r[X]$ be defined by $\mathbb{E}[X(X-1)\ldots (X-r+1)]$ of $X$ for $r\in\mathbb{N}$.

\noindent
The first lemma is based on an idea from \cite{alonspencer}.

\begin{lemma}\label{dom min degree}
Let $G$ be a graph, and let $W$ be a set of vertices of $G$.
If $U$ is a set of vertices such that $|U|$ is minimal, $U\cup W$ is a dominating set of $G$
and $d_G(v)\geq \delta$ for every $v\in V(G)\setminus W$,
then
\begin{align*}
	|U|\leq \frac{1+\ln(\delta+1)}{\delta+1}(n(G)-|W|).
\end{align*}
\end{lemma}

\noindent
\textit{Proof:} 
Let $W\subseteq V(G)$ be given and let $G'=G[V(G)\setminus W]$.
We construct a random set $U$ such that
$U\cup W$ is a dominating set of $G$.
Let $p\in(0,1]$ and 
let every vertex of $G'$ be independently in $U_0$ with probability $p$.
Furthermore, add every vertex of $G'$ in $U_1$ that is not already in $U_0$ and has no neighbor in $U_0\cup W$.
Let $U=U_0\cup U_1$.
Hence 
\begin{align*}
\mathbb{E}[|U|]
=\mathbb{E}[|U_0|]+\mathbb{E}[|U_1|]
&\leq pn(G')+\sum_{v\in V(G')}(1-p)^{d_G(v)+1}\\
&\leq pn(G')+(1-p)^{\delta +1}n(G')\\
&\leq n(G')\left(p+e^{-p(\delta+1)}\right).
\end{align*}
If $p=\frac{\ln (\delta+1)}{\delta+1}$,
then $\mathbb{E}[|U|]\leq \frac{1+\ln(\delta+1)}{\delta+1}n(G')$.
By the first-moment-method,
there is a set $U$ such that $|U|\leq \frac{1+\ln(\delta+1)}{\delta+1}(n(G)-|W|)$ and  $U\cup W$ dominates $G$.
$\Box$

\begin{lemma}\label{jdom}
If $G$ is a graph such that the order of every component is at least $\ell+1$,
then
$$\gamma^{\ell}(G)\leq \frac{n(G)}{\ell+1}.$$
\end{lemma}

\noindent
The proof of Lemma \ref{jdom} is easily done by induction on the number of vertices and thus we omit it.

\begin{thm}[Tillich \cite{til}]\label{thm til}
If $G$ is a connected graph on $k$ vertices,
then there is a constant $c=c(G)>0$ such that
\begin{align*}%\label{boundary}
	b_{G^n}(s)\geq cs(n-\log_k s)
\end{align*}
for all $n\geq 1$ and $1\leq s \leq k^n$.
\end{thm}

\begin{thm}[see for example Durrett \cite{durrett}]\label{thm durret}
Let $X_1,X_2,\ldots$ be a sequence of independent random variables.
If $\mathbb{E}_r[X_n]\rightarrow \lambda^r$ ($n\rightarrow\infty$) for every $r\in \mathbb{N}$ and some $\lambda>0$,
then $X_n$ converges in distribution to a Poisson-distributed random variable with parameter $\lambda$.
\end{thm}

\section{Results and Proofs}\label{section proof}

%%%%%%%%%%%%%%%%%%%%%%%%%%%

The main part of our proof of Theorem \ref{mainresult1} is contained in the following two lemmas.
For this section we fix a sequence $G=(G_n)_{n\geq 1}$ of graphs that satisfies the basic conditions for some $k\in \mathbb{N}\setminus \{1\}$
and fix some sequence $\lambda_n\rightarrow \lambda>0$ as $n\rightarrow\infty$.
We abbreviate $\delta(G_n)$ and $\Delta(G_n)$ by $\delta$ and $\Delta$, respectively. 
Let $p=p(n)$ be chosen such that $\mathbb{E}[X_n]=\lambda$ and let $q=1-p$.
Since $\mathbb{E}[X_n]=\sum_{v\in V(G_n)}q^{d_{G_n}(v)}=\lambda$,
we conclude
\begin{align}\label{magnitude q1}
	q \leq \left(\frac{\lambda_n}{k^n}\right)^{\frac{1}{\Delta}}.
\end{align}
We will frequently use 
\begin{align}\label{magnitude q}
	q^{\delta}
	\leq \left(\frac{\lambda_n}{k^n}\right)^{\frac{1}{n^{1-\epsilon'}}}
	\leq \left( \frac{\lambda_n^{\frac{1}{n}}}{k}\right)^{n^{\epsilon'}}.
\end{align}
Remind that $X_n$ is the (random) number of isolated vertices in the (random) graph $(G_n)_p$.
Note that $\mathbb{E}_r[X_n]$ is the expected number of $r$-tuples of distinct vertices of $G^n_p$
that are all isolated.
In face of our statement we assume that $n$ is sufficiently large 
and note that some inequalities are only true if $n$ is large enough.

\begin{lemma}\label{poi lemma}
Let $G=(G_n)_{n\geq 1}$ be a sequence of connected graphs that satisfies the basic conditions,
$k\in \mathbb{N}\setminus \{1\}$ and let $\lambda_n\rightarrow \lambda>0$ as $n\rightarrow\infty$.
If $p=p(n)$ is chosen such that $\mathbb{E}[X_n]=\lambda$,
then $X_n$ converges in distribution to a Poisson-distributed random variable with parameter $\lambda$.
In particular,
\begin{align*}
	\lim_{n\rightarrow\infty}\mathbb{P}[X_n=0]=e^{-\lambda}.
\end{align*}
\end{lemma}

\noindent
\textit{Proof:}
Let $r\in \mathbb{N}$.
We will show that $\lim_{n\rightarrow\infty}\mathbb{E}_r[X_n]=\lambda^r$.
For $r=1$, there is nothing to show and hence we assume $r\geq 2$.
Let $A_r$ be the set of all $r$-tuples of distinct vertices of $G_n$.
We partition this set into two subsets $B_r$ and $C_r$.
Let $B_r$ be the set of all $r$-tuples such that at least two vertices are adjacent in $G_n$,
that is $B_r=\{(v_1,\ldots,v_r)\in A_r: m(G_n[\{v_1,\ldots,v_r\}])\geq 1\}$ and $C_r=A_r\setminus B_r$.
Note that the vertices in an $r$-tuple of $C_r$ behave independently in view of being isolated vertices in $(G_n)_p$.
Next we show that the $r$-tuples in $B_r$ do not contribute an essential part to the value of $\mathbb{E}_r[X_n]$.
Since $m(G_n[\{v_1,\ldots,v_r\}])\leq \frac{r^2}{2}$, 
we observe the following
\begin{align*}
	\sum_{(v_1,\ldots,v_r)\in B_r}
	\mathbb{P}&[d_{G_n}(v_1)=\ldots=d_{G_n}(v_r)=0]\\
	&= \sum_{(v_1,\ldots,v_r)\in B_r} q^{d_{G_n}(v_1)+\ldots+d_{G_n}(v_r)-m(G_n[\{v_1,\ldots,v_r\}])}\\
	&\leq q^{-\frac{r^2}{2}}\sum_{(v_1,\ldots,v_r)\in B_r} q^{d_{G_n}(v_1)+\ldots+d_{G_n}(v_r)}\\
	&\leq q^{-\frac{r^2}{2}}\sum_{(v_1,\ldots,v_{r-1})\in V(G_n)^{r-1}}\quad
	\sum_{v_r\in N_{G_n}(v_i) {\rm\ for\ some\ } i\in[r-1]} q^{d_{G_n}(v_1)+\ldots+d_{G_n}(v_r)}\\
	&\leq q^{-\frac{r^2}{2}}\sum_{(v_1,\ldots,v_{r-1})\in V(G_n)^{r-1}}\left( q^{d_{G_n}(v_1)+\ldots+d_{G_n}(v_{r-1})} \cdot(r-1)\Delta q^{\delta}\right)\\
	&\leq (r-1)\Delta q^{\delta-\frac{r^2}{2}}\sum_{(v_1,\ldots,v_{r-1})\in V(G_n)^{r-1}} q^{d_{G_n}(v_1)+\ldots+d_{G_n}(v_{r-1})}\\
	&\leq \lambda_n^{r-1}(r-1)\Delta q^{\frac{\delta}{2}}\\
	%&\leq \Delta^2 \left( \frac{\lambda}{k^n}\right)^{ \frac{\delta}{2\Delta}}\\
	&\stackrel{(\ref{magnitude q})}{\leq} \lambda_n^{r-1}(r-1)n^{c}\left( \frac{\lambda}{k^n}\right)^{ \frac{n^{\epsilon'}}{2}}
	=o(1).	
\end{align*}
We establish a lower and an upper bound for the contribution of the elements in $C_r$ to the value of $\mathbb{E}_r[X_n]$.
We have
\begin{align*}
	\sum_{(v_1,\ldots,v_r)\in C_r}
	\mathbb{P}&[d_{G_n}(v_1)=\ldots=d_{G_n}(v_r)=0]\\
	&= \sum_{(v_1,\ldots,v_r)\in C_r} q^{d_{G_n}(v_1)+\ldots+d_{G_n}(v_r)}\\
	&= \sum_{(v_1,\ldots,v_r)\in A_r} q^{d_{G_n}(v_1)+\ldots+d_{G_n}(v_r)}-\sum_{(v_1,\ldots,v_r)\in B_r} q^{d_{G_n}(v_1)+\ldots+d_{G_n}(v_r)}\\
	&\geq \sum_{(v_1,\ldots,v_r)\in V(G_n)^r} q^{d_{G_n}(v_1)+\ldots+d_{G_n}(v_r)}
	- \left((k^{nr}- (k^n)_r) q^{\delta r}\right)
	- (r-1)\Delta q^{\delta}\lambda_n^{r-1}\\
	&\geq \lambda_n^r-r!2^r k^{n(r-1)}\left(\frac{\lambda_n}{k^n}\right)^{\frac{\delta r}{\Delta}}- o(1)\\
	&\geq \lambda_n^r-r!2^r k^{n(r-1)}\left(\frac{\lambda_n}{k^n}\right)^{rn^{\epsilon'}}- o(1)\\
	&= \lambda^r_n - o(1).
\end{align*}
Furthermore,
\begin{align*}
  \sum_{(v_1,\ldots,v_r)\in C_r}
	\mathbb{P}[d_{G_n}(v_1)=\ldots=d_{G_n}(v_r)=0]
	&= \sum_{(v_1,\ldots,v_r)\in C_r} q^{d_{G_n}(v_1)+\ldots+d_{G_n}(v_r)}\\
	&\leq \sum_{(v_1,\ldots,v_r)\in V(G_n)^r} q^{d_{G_n}(v_1)+\ldots+d_{G_n}(v_r)}\\
	&= \lambda_n^r. 
\end{align*}
Combining these three bounds,
we conclude
\begin{align*}
	\lambda^r_n - o(1)
	\leq \mathbb{E}_r[X_n]
	\leq \lambda^r_n +o(1)
\end{align*}
and hence $\lim_{n\rightarrow\infty}\mathbb{E}_r[X_n]=\lambda^r$.
By Theorem \ref{thm durret}, this completes the proof of Lemma \ref{poi lemma}.
$\Box$

%%%%%%%%%%%%%%%%%%%%%%%%%%%%%%%%%%%%%%%%%%%%
%%%%%%%%%%%%%%%%%%%%%%%%%%%%%%%%%%%%%%%%%%%%

\begin{lemma}\label{main lemma}
Let $G=(G_n)_{n\geq 1}$ be a sequence of connected graphs that satisfies the basic conditions,
$k\in \mathbb{N}\setminus \{1\}$ and let $\lambda_n\rightarrow \lambda>0$ as $n\rightarrow\infty$.
If $p=p(n)$ is chosen such that $\mathbb{E}[X_n]=\lambda$,
then
\begin{align*}
	\lim_{n\rightarrow\infty}\mathbb{P}\left[(G_n)_p {\rm\ has\ a\ component\ of\ order\ } 2\leq s\leq \frac{k^n}{2}\right]=0.
\end{align*}
\end{lemma}

\noindent
\textit{Proof:}
We frequently use the well known inequality $\binom{n}{k}\leq \left(\frac{en}{k}\right)^k$ and
omit roundings for more clarity.
The reader may convince himself that all inequalities are still correct if we add all necessary roundings.
Throughout the proof we denote by $S$ a set of vertices of $G_n$ that may form a component in $G_p^n$ and by $s$ its cardinality. 
%Let $\mathcal{A}_s$ be the set of all subsets of vertices of $G^n_p$ such that $|S|=s$,
%that is, $\mathcal{A}_s =\{S\subseteq V((G_n)):|S|=s\}$.
One key part of the proof is to find a good upper bound for the number of connected components in $G_n$ of order $s$.
The proof is divided into a three cases.
In the first two cases we consider small values of $s$ and in the third we assume that $s$ is large.
We do not start with an upper bound of the number of connected components in $G_n$ of order $s$,
but with an upper bound for the number of connected components in $G_n$ of order $s$ containing some vertex $v\in V(G_n)$.
Let $S$ be the set of vertices of such a component.
Since $v$ is in $S$ and $S$ is connected,
there is an ordering $v_1v_2\ldots v_s$ of $S$
such that $v=v_1$ and $v_i$ is adjacent to some $v_j$ and $j<i$ for all $i\in\{2,\ldots s\}$.
Thus
\begin{align*}
	\left|\left\{ S\in\binom{V(G_n)}{s}: G_n[S] {\rm \ is\ connected\ and\ }v\in S\right\} \right| \\
	< \Delta \cdot (2\Delta ) \cdot  \ldots \cdot ((s-1)\Delta)
	\leq (\Delta s)^s.
\end{align*}
%By Theorem \ref{thm til}, there is a positive constant $z=z(G)$,
%which depends only on $G$,
%such that $b_{G_n}(s)\geq zs(n-\log_k s)$
%for all $n\geq 1$ and $1\leq s \leq k^n$.

\noindent
\textbf{Case 1} $(2\leq s \leq n^{1- \epsilon'})$:

\noindent
Note that every vertex in $G_n$ has degree at least $\delta$ and
$m(G_n[S])\leq \frac{s^2}{2}$.
If $v\in S$,
then $b_{G_n}(S)\geq d_{G_n}(v)+ \delta(s-1)- s^2\geq d_{G_n}(v)+ \frac{\delta(s-1)}{2}$.
Hence

\begin{align*}
	\sum_{S\in\binom{V(G_n)}{s}} \mathbb{P}[(G_n)_p[S] {\rm\ is\ a\ component\ in\ }(G_n)_p]
	&\leq \sum_{v\in V(G_n)}(\Delta s)^s q^{d_{G_n}(v)+ \frac{\delta(s-1)}{2}}\\
	&= \lambda_n(\Delta s)^s q^{\frac{\delta (s-1)}{2}}\\
	%&\leq (\Delta s)^s \left(\frac{\lambda^{\frac{1}{n}}}{k}\right)^{\frac{n}{\Delta}\frac{\delta(s-1)}{2}}\\
	&\stackrel{(\ref{magnitude q})}{\leq} \lambda_n\Delta s\left(\Delta s\left(\frac{\lambda_n^{\frac{1}{n}}}{k}\right)^{\frac{n^{\epsilon'}}{2}}\right)^{s-1}\\
	&\leq \lambda_n n^{c+1}\left(\frac{1}{n^{3c}}\right)^{s-1}\\
	&\leq \left(\frac{1}{n^{c}}\right)^{s-1}.
\end{align*}
Thus\begin{align*}
	\sum_{s=2}^{n^{1- \epsilon'}} \sum_{S\in\binom{V(G_n)}{s}} \mathbb{P}[(G_n)_p[S] {\rm\ is\ a\ component\ in\ }(G_n)_p]
	=o(1)\quad (n\rightarrow\infty).
\end{align*}

\noindent
\textbf{Case 2} $(n^{1- \epsilon'}+1\leq s \leq k^{\frac{\epsilon}{4}n})$:

\noindent
%Note that $z\leq\Delta$ by plugging $s=1$ into inequality (\ref{boundary}).
%Thus $s\leq k^{\frac{1}{4}n}$.
%Let $r=r(G,z)=\frac{4\Delta}{z}$ 
Since $s\leq k^{\frac{\epsilon}{4}n}$,
we obtain $\epsilon\Delta s\left(1- \frac{1}{n}\log_k s\right)-\Delta\geq \frac{\epsilon\Delta s}{2}$.
Using a similar idea as in Case 1,
we obtain
\begin{align*}
	\sum_{S\in\binom{V(G_n)}{s}} \mathbb{P}[(G_n)_p[S] {\rm\ is\ a\ component\ in\ }(G_n)_p]
	&\leq \sum_{v\in V(G_n)}(\Delta s)^s q^{b_{G_n}(s)}\\
	&\leq \sum_{v\in V(G_n)}(\Delta s)^s q^{d_{G_n}(v)-\Delta +b_{G_n}(s)}\\
	&\leq \lambda_n(\Delta s)^s q^{\epsilon\Delta s\left(1-\frac{1}{n}\log_ks\right)-\Delta}\\
	&\leq \lambda_n(\Delta s)^s q^{\frac{\epsilon\Delta s}{2}}.
\end{align*}
Let $f(s)=\lambda_n(\Delta s)^s q^{\frac{\epsilon \Delta s}{2}}$ and hence
$f'(s)=f(s)\ln\left(e\Delta s q^{\frac{\epsilon\Delta}{2}}  \right)$.
If $e\Delta s q^{\frac{\epsilon\Delta}{2}}  <1$,
then $f(s)$ is monotone decreasing in $s$.
Since $s\leq k^{\frac{\epsilon}{4}n}$,
we have 
\begin{align*}
	e\Delta s q^{\frac{\epsilon\Delta}{2}}
	<e\Delta  k^{\frac{\epsilon}{4}n}q^{\frac{\epsilon\Delta}{2}}
	\stackrel{(\ref{magnitude q})}{\leq} e\Delta \lambda^{\frac{\epsilon}{4}} q^{\frac{\epsilon\Delta}{4}}
	\leq 1.
\end{align*}
Therefore,
\begin{align*}
	\sum_{s=n^{1- \epsilon'}+1}^{k^{\frac{\epsilon}{4}n}} \sum_{S\in\binom{V(G_n)}{s}} \mathbb{P}[(G_n)_p[S] {\rm\ is\ a\ component\ in\ }(G_n)_p]
	&\leq k^{\frac{\epsilon}{4}n} \lambda_n(\Delta n^{1- \epsilon'})^{n^{1- \epsilon'}} q^{\frac{\epsilon\Delta n^{1- \epsilon'}}{2}}\\
	&\stackrel{(\ref{magnitude q1})}{\leq} \lambda_n n^{(c+1)n} \left(\frac{\lambda_n}{k^n}\right)^{\frac{\epsilon n^{1- \epsilon'}}{2}}\\
	&=o(1)\quad (n\rightarrow\infty).
\end{align*}

\noindent
\textbf{Case 3} $(k^{\frac{\epsilon}{4}n}\leq s \leq k^{n}/2)$:

\noindent
We partition $\binom{V(G_n)}{s}$ into two sets $\mathcal{A}_s$ and $\mathcal{B}_s$,
where 
$$\mathcal{A}_s=\left\{S\in \binom{V(G_n)}{s}: b_{G_n}(S)\geq \frac{\epsilon \Delta s}{2}\left(1-\frac{1}{n}\log_k s +\frac{1}{n}\ln^2 n\right) \right\}$$ 
and
$\mathcal{B}_s=\binom{V(G_n)}{s}\setminus \mathcal{A}_s$.
First, we only consider $S\in \mathcal{A}_s$.
To establish a new upper bound for the number of connected components $S$ of $G_n$
we argue as follows.
By Lemma \ref{jdom}, there is a subset $U$ of $\frac{\epsilon}{2}s$ vertices of $S$
that $(\frac{2}{\epsilon}-1)$-dominates $S$.
There are at most $\binom{k^n}{\frac{\epsilon}{2}s}$ choices for such a set.
The remaining vertices of $S$ are located close to at least one vertex in $U$.
In detail, for every vertex in $S$ there is a vertex in $U$ in distance at most $\frac{2}{\epsilon}-1$.
Note that for every vertex $v$ of $G_n$ there are at most $\Delta^{\frac{2}{\epsilon}}$ vertices in distance at most $\frac{2}{\epsilon}-1$ to $v$.
This leads to
\begin{align*}
	\left|\left\{S\in\binom{V(G_n)}{s}: G_n[S] {\rm\ is\ connected}\right\}\right|
	\leq \binom{k^n}{\frac{\epsilon}{2}s}\binom{\frac{\epsilon}{2}\Delta^{\frac{2}{\epsilon}}s}{\frac{2- \epsilon}{2}s}
\end{align*}
and hence
\begin{align*}
	\sum_{S\in\mathcal{A}_s} \mathbb{P}[(G_n)_p[S] {\rm\ is\ a\ component\ in\ }&(G_n)_p]\\
	&\leq \binom{k^n}{\frac{\epsilon}{2}s}\binom{\frac{\epsilon}{2}\Delta^{\frac{2}{\epsilon}}s}{\frac{2- \epsilon}{2}s}
	q^{\frac{\epsilon \Delta s}{2}\left(1-\frac{1}{n}\log_k s +\frac{1}{n}\ln^2 n\right)}\\
	&\stackrel{(\ref{magnitude q1})}{\leq} \left( \frac{2ek^n}{\epsilon s} \right)^{\frac{\epsilon}{2}s}
	\left( \frac{\epsilon e\Delta^{\frac{2}{\epsilon}}}{2- \epsilon} \right)^{\frac{2-\epsilon}{2}s}
	\left( \frac{\lambda_n}{k^n} \right)^{\frac{\epsilon s}{2}\left(1-\frac{1}{n}\log_k s +\frac{1}{n}\ln^2 n\right)}\\
	&\leq \left( \frac{2ek^n}{\epsilon s} \right)^{\frac{\epsilon}{2}s}
	\left( \frac{(\epsilon e)^{\frac{2- \epsilon}{\epsilon}}\Delta^{\frac{2(2- \epsilon)}{\epsilon^2}}}{(2- \epsilon)^{\frac{2- \epsilon}{\epsilon}}} \right)^{\frac{\epsilon}{2}s}
	\left( \frac{2\lambda_n s}{k^n n^{\ln k \ln n}} \right)^{\frac{\epsilon }{2}s}\\
	&\leq \left( \frac{1}{n} \right)^{s}.
\end{align*}
This implies that
\begin{align*}
	\sum_{s=k^{\frac{\epsilon}{4}n}}^{\frac{k^n}{2}} \sum_{S\in\mathcal{A}_s} \mathbb{P}[(G_n)_p[S] {\rm\ is\ a\ component\ in\ }(G_n)_p]
	&\leq \sum_{s=k^{\frac{\epsilon}{4}n}}^{\frac{k^n}{2}}\left( \frac{1}{n} \right)^{s} 
	=o(1)\quad (n\rightarrow\infty).
\end{align*}
If $b_{G_n}(s)\geq \frac{\epsilon \Delta s}{2}\left(1-\frac{1}{n}\log_k s +\frac{1}{n}\ln^2 n\right)$,
then $\mathcal{B}_s=\emptyset$.
%Thus we only consider these values for $s$ if $b_{G_n}(s)\leq \frac{c_1 s}{j}\left(n-\log_k \frac{s}{n^j}\right)$.
Since
\begin{eqnarray*}
	& \epsilon\Delta s\left(1-\frac{1}{n}\log_k s \right) \leq b_{G_n}(s) &\leq \frac{\epsilon \Delta s}{2}\left(1-\frac{1}{n}\log_k s +\frac{1}{n}\ln^2 n\right)\\
	\Rightarrow& n-\log_k s  &\leq \ln^2 n\\
	%\Rightarrow& \frac{k^n}{s} &\leq n^{\ln k \ln n}\\
	\Rightarrow& \frac{k^n}{n^{\ln k \ln n}}& \leq s,
\end{eqnarray*}
from now on we may assume that $s\geq \frac{k^n}{n^{\ln k \ln n}}$ and $S\in \mathcal{B}_s$ and hence
\begin{align*}
b_{G_n}(S)
\leq \frac{\epsilon \Delta s}{2}\left(1-\frac{1}{n}\log_k s +\frac{1}{n}\ln^2 n\right)
\leq \frac{\epsilon\Delta  \ln^2 n}{n}s.
\end{align*}
Let $W$ be a set of those vertices in $S$ that have
at least $\delta n^{- \frac{\epsilon}{2}}$ neighbors in $G_n[V(G_n)\setminus S]$ and let $w=|W|$.
By double counting the edges of the boundary of $S$, we obtain
\begin{align*}
w\leq \frac{\epsilon \Delta  \ln^2 n}{\delta n^{1- \frac{\epsilon}{2}}}s
\leq  \epsilon  n^{- \frac{\epsilon}{2}}\ln^2 ns.
\end{align*}
By the choice of $W$, 
every vertex in $S\setminus W$ has degree at least $\delta \left(1-n^{- \frac{\epsilon}{2}}\right)$ in $G_n[S]$.
Let a set of vertices $U$ together with $W$ be a minimal dominating set of $G_n[S]$ under the condition that $W$ is already given
and let $u=|U|$.
By Lemma \ref{dom min degree}, we conclude
\begin{align*}
u\leq \frac{1+ \ln\left(\delta \left(1-n^{- \frac{\epsilon}{2}}\right)\right)}{\delta \left(1-n^{- \frac{\epsilon}{2}}\right)}s
\leq \frac{2\ln \delta}{\delta}s.
\end{align*}
Note that every vertex of $S$ is in $U$, in $W$ or in the neighborhood of $U$.
Furthermore, every vertex in $U$ has at most $\delta n^{- \frac{\epsilon}{2}}$ neighbors
not in $S$.
Thus, to get an upper bound on $|\mathcal{B}_s|$,
we first choose the vertices of $U$ and $W$ arbitrarily and then the non-neighbors of $U$.
Hence,
\begin{align*}
|\mathcal{B}_s|
&\leq \binom{k^n}{u}
\binom{k^n}{w }
\left(\sum_{(k_1,\ldots,k_u)\in\{0,\ldots,\delta n^{- \frac{\epsilon}{2}}\}^u}\prod_{i=1}^{u}\binom{\Delta}{k_i} \right)\\
&\leq \left(\frac{e \delta k^n}{2 \ln \delta s}\right)^{\frac{2 \ln \delta }{\delta }s}
\left(\frac{en^{\frac{\epsilon}{2}}k^n}{\epsilon \ln^2 n s}\right)^{\frac{\epsilon \ln^2 n}{n^{\frac{\epsilon}{2}}}s}
\left(\delta n^{- \frac{\epsilon}{2}}+1\right)^u 
{\binom{\Delta}{\delta n^{- \frac{\epsilon}{2}}}}^u\\
&\leq \left(n^{2\ln k \ln n}\right)^{\frac{2 \ln \delta }{\delta }s}
\left(n^{2\ln k \ln n}\right)^{\frac{\epsilon \ln^2 n}{n^{\frac{\epsilon}{2}}}s}
\cdot \left(n^{c+1}\right)^{\frac{2 \ln \delta }{\delta }s}
\cdot n^{\frac{2\ln \delta}{ n^{\frac{\epsilon}{2}}}s}\\
&\leq n^{\frac{1}{n^{\frac{\epsilon}{4}}}s}\\
&=\left(\exp\left(\frac{1}{n^{\frac{\epsilon}{4}}}\ln n \right)\right)^s\\
&\leq (1+\epsilon_1)^s
\end{align*}
for every fixed $\epsilon_1>0$.
We choose $\epsilon_1$ small enough such that there is an $\epsilon_2>0$
such that $(1+\epsilon_1)\left(\frac{\lambda_n^{\frac{1}{n}}}{k}\right)^{\epsilon\log_k 2}\leq 1- \epsilon_2$.
Since $s\leq \frac{k^n}{2}$, 
we have $b_{G_n}(s)\geq \frac{\epsilon\Delta s \log_k 2}{n}$.
Therefore, 
\begin{align*}
	\sum_{s=\frac{k^n}{n^{\ln k \ln n}}}^{\frac{k^n}{2}} \sum_{S\in\mathcal{B}_s} \mathbb{P}[(G_n)_p[S] {\rm\ is\ a\ component\ in\ }(G_n)_p]
	&\leq \sum_{s=\frac{k^n}{n^{\ln k \ln n}}}^{\frac{k^n}{2}}|\mathcal{B}_s|q^{b_{G_n}(s)}\\
	&\stackrel{(\ref{magnitude q1})}{\leq} \sum_{s=\frac{k^n}{n^{\ln k \ln n}}}^{\frac{k^n}{2}}(1+\epsilon_1)^s \left(\left(\frac{\lambda_n^{\frac{1}{n}}}{k}\right)^{\epsilon\log_k 2}\right)^s\\
	&\leq \sum_{s=\frac{k^n}{n^{\ln k \ln n}}}^{\frac{k^n}{2}}(1-\epsilon_2)^s\\
	&=o(1)\quad (n\rightarrow \infty).
\end{align*}
This completes the proof of Lemma \ref{main lemma}.
$\Box$

\bigskip

\noindent
\textit{Proof of Theorem \ref{mainresult1}:} 
Note that every disconnected graph without isolated vertices has a component of
order between 2 and half of the order of the graph.
By Lemma \ref{main lemma}, we obtain
\begin{align*}
  0\leq& \lim_{n\rightarrow\infty}\mathbb{P}[(G_n)_p {\rm\ is\ disconnected}]-\mathbb{P}[X_n>0]\\
  \leq&\lim_{n\rightarrow\infty} \mathbb{P}\left[(G_n)_p {\rm\ has\ a\ component\ of\ order\ } 2\leq s\leq \frac{k^n}{2}\right]
  =0,
\end{align*}
and hence
\begin{align*}
  \lim_{n\rightarrow\infty}\mathbb{P}[(G_n)_p {\rm\ is\ disconnected}]=\lim_{n\rightarrow\infty}\mathbb{P}[X_n>0].
\end{align*}
By Lemma \ref{poi lemma}, we conclude
\begin{align*}
  \lim_{n\rightarrow\infty}\mathbb{P}[X_n>0]=1-e^{-\lambda}.
\end{align*}
$\Box$

\bigskip

\noindent
\textit{Proof of Corollary \ref{main result}:}
We first verify that the sequence $(H^n)_{n\geq 1}$ of Cartesian powers of $H$ satisfies the basic conditions.
Let $k=n(H^1)$ and hence condition 1 is trivially satisfied.
For a vertex $v=(v_1,\ldots,v_n)\in V(H^n)$, 
we have $d_{H^n}(v)=d_H(v_1)+\ldots+d_H(v_n)$.
Thus $\delta(H^n)=n\delta(H)$ and $\Delta(H^n)=n\Delta(H)$ and 
hence $\frac{\Delta(H^n)}{\delta(H^n)}=\frac{\Delta(H)}{\delta(H)}$.
Thus the Conditions 2, 3 and 4 are satisfied and by Theorem \ref{thm til} Condition 5 is satisfied.
Recall that $\alpha_{i,H}$ is the number of vertices of degree $i$ in $H$.
Furthermore,
\begin{align*}
	\lambda_n=\mathbb{E}[X_n]
	&=\sum_{v\in V(H^n)}q^{d_{H^n}(v)}\\
	&=\sum_{(v_1,\ldots,v_n)\in V(H^n)}q^{d_{H^n}(v_1)+\ldots+d_{H^n}(v_n)}\\
	&=\left(\sum_{u\in V(H)}q^{d_{H}(u)}\right)^n\\
	&=\left(\sum_{i=1}^{\infty}\alpha_{i,H} q^i\right)^n.
\end{align*}
Thus $q$ is the unique positive solution of the equation $P(H,x)=\sum_{i=1}^{\infty}\alpha_{i,H} q^i =\lambda_n^{\frac{1}{n}}$.
$\Box$

\section{Conclusion}
Theorem \ref{mainresult1} implies that the threshold function for sparse graphs do only depend on the degree sequence of the graph
if there is no small set of edges that disconnects the graph.
We conjecture that a more general versions of Theorem \ref{mainresult1} is true,
that is, that the Conditions 3, 4 and 5 can be weakened.
Maybe other techniques are necessary to prove such a result. 

\section{Acknowledgment}
The author thanks Zakhar Kabluchko for introducing him into this topic and
Dieter Rautenbach for valuable discussions.


\begin{thebibliography}{1}

\bibitem{alonspencer}
N. Alon and J.H. Spencer, The probabilistic method, Third Edition, Wiley-Interscience Series in Discrete Mathematics and Optimization,
John Wiley \& Sons Inc, Hoboken, NY, 2008.

\bibitem{bolobook}
B. Bollob{\'a}s. Random graphs. Second edition. Cambridge Studies in Advanced Mathematics, 73. Cambridge University Press, Cambridge, 2001.

\bibitem{bolo}
B. Bollob{\'a}s, The evolution of the cube, In \textit{Combinatorial mathematics ({M}arseille-{L}uminy, 1981)},
volume 75 of \textit{North-Holland Math. Stud.}, pages 91-97, North-Holland, Amsterdam, 1983.

\bibitem{burtin}
J.D. Burtin, The probability of connectedness of a random subgraph of an {$n$}-dimensional cube, Problemy Pereda\v ci Informacii \textbf{13} (1977) 90-95.

\bibitem{clark}
L. Clark, Random subgraphs of certain graph powers, Int. J. Math. Math. Sci. \textbf{32} (2002) 285-292.

\bibitem{durrett}
R. Durrett. \textit{Probability}. The Wadsworth \& Brooks/Cole Statistics/Probability Series.
Wadsworth \& Brooks/Cole Advanced Books \& Software, Pacific Grove, CA, 1991.
Theory and examples.

\bibitem{erdosrenyi}
P. Erd{\H{o}}s, and A. R{\'e}nyi, On random graphs. {I}, Publ. Math. Debrecen \textbf{6} (1959) 290-297.

\bibitem{erdosspencer}
P. Erd{\H{o}}s and J.H. Spencer, Evolution of the {$n$}-cube, Comput. Math. Appl. \textbf{5} (1979) 33-39.

\bibitem{jo}
F. Joos, Random subgraphs in Cartesian powers of regular graphs, Electron. J. Combin. \textbf{9} (2012) Paper 47.

\bibitem{hofstad}
R. van der Hofstad and A. Nachmias, Hypercube percolation, preprint, arXiv:1201. 3953v1.

\bibitem{palasti}
I. Pal{\'a}sti, On the connectedness of random graphs, In \textit{Studies in {M}athematical {S}tatistics: {T}heory and {A}pplications},
pages 105-108, Akad. Kiad\'o, Budapest 1968.

\bibitem{ruci}
A. Ruci{\'n}ski, The {$r$}-connectedness of {$k$}-partite random graph, Bull. Acad. Polon. Sci. S\'er. Sci. Math. \textbf{29} (1981) 321-330.

\bibitem{til}
J.-P. Tillich, Edge isoperimetric inequalities for product graphs, Discrete Math. \textbf{213}(1-3) (2000) 291-320.

\end{thebibliography}
\end{document}